\documentclass[12pt]{amsart}
\usepackage{amsmath,amscd,amssymb,amsfonts}
\newcommand{\h}{\hbox}
\newcommand{\q}{\quad}
\newcommand{\nin}{\par\noindent}
\newcommand{\bs}{\par\bigskip}
\newcommand{\ms}{\par\medskip}
\newcommand{\sk}{\par\smallskip}
\newcommand{\bsn}{\par\bigskip\noindent}
\newcommand{\msn}{\par\medskip\noindent}

\newcommand{\bl}{\bigl}
\newcommand{\br}{\bigl}
\newcommand{\msum}{\h{$\sum$}}
\newcommand{\mopl}{\h{$\bigoplus$}}
\newcommand{\ssb}{\raise.15ex\h{${\scriptscriptstyle\bullet}$}}
\newcommand{\A}{{\mathcal A}}
\newcommand{\al}{\alpha}
\newcommand{\af}{\mathfrak{a}}
\newcommand{\alt}{\widetilde{\alpha}}
\newcommand{\C}{{\mathbb C}}
\newcommand{\dd}{\partial}
\newcommand{\e}{{\mathbf e}}
\newcommand{\B}{{\mathcal B}}
\newcommand{\BB}{\widetilde{\mathcal B}}
\newcommand{\D}{{\mathcal D}}
\newcommand{\ep}{\varepsilon}
\newcommand{\G}{{\mathcal G}}
\newcommand{\Hc}{{\mathcal H}}
\newcommand{\Ht}{\widetilde{H}}
\newcommand{\I}{{\mathcal I}}
\newcommand{\J}{{\mathcal J}}
\newcommand{\JC}{{\rm JC}}
\newcommand{\JCt}{\widetilde{\rm JC}}
\newcommand{\Jt}{\widetilde{\mathcal J}}
\newcommand{\Gt}{\widetilde{\mathcal G}}
\newcommand{\M}{\widetilde{M}}
\newcommand{\la}{\lambda}
\newcommand{\lct}{{\rm lct}}
\newcommand{\N}{{\mathbb N}}
\newcommand{\OO}{{\mathcal O}}
\newcommand{\Q}{{\mathbb Q}}
\newcommand{\R}{{\mathbb R}}
\newcommand{\Si}{\Sigma}
\newcommand{\Vt}{\widetilde{V}}
\newcommand{\X}{\widetilde{X}}
\newcommand{\om}{\omega}
\newcommand{\omt}{\widetilde{\omega}}
\newcommand{\Z}{{\mathbb Z}}
\newcommand{\can}{{\rm can}}
\newcommand{\DR}{{\rm DR}}
\newcommand{\Gr}{{\rm Gr}}
\newcommand{\Sing}{{\rm Sing}}
\newcommand{\into}{\hookrightarrow}
\newcommand{\simto}{\buildrel\sim\over\longrightarrow}
\newcommand{\ges}{\geqslant}
\newcommand{\les}{\leqslant}
\newcommand{\1}{\hskip1pt}
\newcommand{\sot}{\1{\otimes}\1}
\begin{document}
\title[Thom-Sebastiani theorems]
{Thom-Sebastiani theorems for filtered $\D$-modules\\and for multiplier ideals}
\author[L. Maxim ]{Laurentiu Maxim}
\address{L. Maxim : Department of Mathematics, University of Wisconsin-Madison, 480 Lincoln Drive, Madison WI 53706-1388 USA}
\email{maxim@math.wisc.edu}
\author[M. Saito ]{Morihiko Saito}
\address{M. Saito: RIMS Kyoto University, Kyoto 606-8502 Japan}
\email{msaito@kurims.kyoto-u.ac.jp}
\author[J. Sch\"urmann ]{J\"org Sch\"urmann}
\address{J. Sch\"urmann : Mathematische Institut, Universit\"at M\"unster, Einsteinstr. 62, 48149 M\"unster, Germany}
\email{jschuerm@uni-muenster.de}
\begin{abstract} We give a proof of the Thom-Sebastiani type theorem for holonomic filtered $D$-modules satisfying certain good conditions (including Hodge modules) by  using algebraic partial microlocalization. By a well-known relation between multiplier ideals and $V$-filtrations of Kashiwara and Malgrange, the argument in the proof implies also a Thom-Sebastiani type theorem for multiplier ideals, which cannot be deduced from a already known proof of the Thom-Sebastiani theorem for mixed Hodge modules (since the latter gives only the information of graded pieces of multiplier ideals). We also sketch a more elementary proof of the Thom-Sebastiani type theorem for multiplier ideals (as communicated to us by M.~Musta\c{t}\v{a}), which seems to be known to specialists, although it does not seem to be stated explicitly in the literature.
\end{abstract}
\maketitle
\centerline{\bf Introduction}
\bsn
For $a=1,2$, let $Y_a$ be a smooth complex variety or a connected complex manifold. Set $X_a:=f_a^{-1}(0)$ with $f_a$ a non-constant function on $Y_a$ (that is, $f_a\in\Gamma(Y_a,\OO_{Y_a})\setminus\C$). Put
$$Y:=Y_1\times Y_2,\q X:=f^{-1}(0)\subset Y\q\h{with}\q f:=f_1+f_2\,\,\,\h{on}\,\,\,Y,$$
where $f_1+f_2$ is an abbreviation of ${\rm pr}_1^*f_1+{\rm pr}_2^*f_2$ with ${\rm pr}_a$ the $a$-th projection $(a=1,2)$. For $\al\in\Q$, we have a decreasing sequence of the {\it multiplier ideals} $\J(\al X)\subset\OO_Y$ together with their graded quotients $\G(\al X)$, and similarly for $\J(\al X_a)$, $\G(\al X_a)$ with $a=1,2$, see also (2.1) below. (For general references to multiplier ideals, see \cite{La}, \cite{Na}.) Set $\Si_a:={\rm Sing}\,X_a$, $\Si:={\rm Sing}\,X$. In this paper we prove the following.
\msn
{\bf Theorem~1.} {\it In the above notation, we have the following equalities for $\al\in(0,1):$
$$\J(\al X)=\msum_{\al_1+\al_2=\al}\,\J(\al_1X_1)\boxtimes\J(\al_2X_2)\q\h{in}\,\,\,\,\OO_Y=\OO_{Y_1}\boxtimes\OO_{Y_2},
\leqno(1)$$
together with the canonical isomorphisms of $\OO_X$-modules for $\al\in(0,1):$
$$\G(\al X)=\msum_{\al_1+\al_2=\al}\,\G(\al_1X_1)\boxtimes\G(\al_2X_2),
\leqno(2)$$
by replacing {$($if necessary$)$} $Y_a$ with an open neighborhood of $X_a$ in $Y_a$ {$(a=1,2)$} so that $\Si=\Si_1\times\Si_2$.}
\ms
Here we may assume $\al_1,\al_2\in(0,\al)$ on the right-hand side of (1), (2).
The formula (1) determines $\J(\al X)$ for any $\al\in\Q$, since it is well known (and is easy to show) that
$$\J((\al+1)X)=f\J(\al X)\,\,\,\,(\forall\,\al\ges 0),\q\J(\al X)=\OO_Y\,\,\,\,(\forall\,\al\les 0),
\leqno(3)$$
see (2.1.2) below. The formula (2) for $\al=1$ is more complicated (see Corollary~(2.4) below), since it is closely related to the ``irrationality" of the singularities of $X$, see a remark after (2.1.9) below.
\sk
We define the set of {\it jumping coefficients} $\JC(X)$ and the {\it log canonical threshold\,} $\lct(X)$ by
$$\JC(X):=\{\al\in\Q\mid\G(\al X)\ne 0\},\q\q\lct(X):=\min\JC(X),$$
and similarly for $\JC(X_a)$, $\lct(X_a)$ ($a=1,2$). From Theorem~1 we can deduce the following {\it addition theorem for jumping coefficients and log canonical thresholds}.
\msn
{\bf Corollary~1.} {\it In the above notation, we have the equalities
$$\JC(X)\cap(0,1)=\bl(\JC(X_1)+\JC(X_2){\br)}\cap(0,1),
\leqno(4)$$
\vskip-6mm
$$\lct(X)=\min\bl\{1,\lct(X_1)+\lct(X_2)\br\},
\leqno(5)$$
by replacing {$($if necessary$)$} $Y_a$ with an open neighborhood of $X_a$ {$(a=1,2)$} so that $\Si=\Si_1\times\Si_2$.}
\ms
Note that $\JC(X)$ is determined by $\JC(X_1)$, $\JC(X_2)$ using (4) together with the periodicity coming from (3) (see also (2.1.5) below), since we have $1\in\JC(X)$ by looking at the smooth points of $X$. Note also that Theorem~1 and Corollary~1 were {\it essentially} known to specialists according to M.~Musta\c{t}\v{a} although it does not seem to be stated explicitly in the literature, see Remark~(2.7) below.
\sk
Theorem~1 is obtained as a by-product of a proof of the Thom-Sebastiani theorem for filtered holonomic $\D$-modules in the constant coefficient case using the {\it algebraic partial microlocalization}, where a well-known relation between the multiplier ideals and the $V$-filtration of Kashiwara \cite{Ka2} and Malgrange \cite{Ma2} is also used, see \cite[Theorem 0.1]{BS} (and also (2.1.6) below). Note also that the equality (1) in Theorem~1 cannot be deduced from the arguments in \cite{ts} (which gives only the information of graded pieces of multiplier ideals). Let $i_f:X\into Y$ be the graph embedding by $f$. Set
$$(\B_f,F):=(i_f)_*^{\!\D}(\OO_X,F),$$
where $(i_f)_*^{\!\D}$ is the direct image as filtered $\D$-module. In this paper, the Hodge filtration $F$ is {\it indexed like a right $\D$-module}. This means that the Hodge filtrations are {\it not shifted\1} by codimensions under the direct images by closed embeddings, see for instance \cite[Section 1.2]{MSS}, \cite[Section B.3]{ypg}. We have the filtration $V$ of Kashiwara \cite{Ka2} and Malgrange \cite{Ma2} on $\B_f$ indexed decreasingly by $\Q$ so that $\dd_tt-\al$ is nilpotent on $\Gr_V^{\al}\B_f$, see \cite{mhp}. Setting $\e(-\al):=\exp(-2\pi i\al)$, we define the $\e(-\al)$-eigenspace of the {\it filtered vanishing cycle functor} by
$$\varphi^{(\al)}_f(\OO_Y,F):=\Gr_V^{\al}(\B_f,F)\q\q\bl(\al\in(-1,0]\br),$$
so that $\mopl_{\al\in(-1,0]}\,\varphi^{(\al)}_f(\OO_Y,F)$ is the underlying $F$-filtered $\D_Y$-module of the mixed Hodge module of vanishing cycles in \cite{mhp}, \cite{mhm} (which is denoted by $\varphi_f\Q_{h,Y}[d_Y-1]$ in this paper), see (1.1.10) below. There are canonical isomorphisms
$$\DR(\varphi^{(\al)}_f\OO_Y)=\varphi_{f,\e(-\al)}\C_Y[d_Y-1],$$
where $d_Y:=\dim Y$, and $\varphi_{f,\e(-\al)}$ denotes the $\e(-\al)$-eigenspace of the vanishing cycle functor for $\C$-complexes \cite{De}, see (1.5.2) below. (Similarly for $\varphi^{(\al)}_{f_a}(\OO_{Y_a},F)$ with $a=1,2$.)
\sk
In this paper we give a proof of the Thom-Sebastiani theorem for holonomic filtered $\D$-modules satisfying the assumptions (1.1.5--7) below (which hold in the case of Hodge modules), see Theorem~(1.2) below. In a special case this implies the following.
\msn
{\bf Theorem~2.} {\it In the above notation, there are canonical isomorphisms of filtered $\D_Y$-modules for $\al\in(-1,0]:$
$$\aligned\varphi_f^{(\al)}(\OO_Y,F)&=\mopl_{\al_1\ges\al}\varphi^{(\al_1)}_{f_1}(\OO_{Y_1},F)\boxtimes\varphi^{(\al-\al_1)}_{f_2}(\OO_{Y_2},F)\\
&\q\oplus\,\mopl_{\al_1<\al}\varphi^{(\al_1)}_{f_1}(\OO_{Y_1},F)\boxtimes\varphi^{(\al-1-\al_1)}_{f_2}(\OO_{Y_2},F[-1]),\endaligned
\leqno(6)$$
by replacing {$($if necessary$)$} $Y_a$ with an open neighborhood of $X_a$ {$(a=1,2)$} so that $\Si=\Si_1\times\Si_2$.}
\ms
More precisely, $\al_1$ on the right-hand side of $(6)$ is contained in $(-1,0]\cap[\al,\al+1)$ or in $(-1,0]\cap[\al-1,\al)$, see (1.2.2) below. (Note that Theorem~2 also follows from \cite{ts}.)
\sk
For the proof of Theorem~2, we use the {\it algebraic partial microlocalization} $\BB_{f_a}$, $\BB_f$ of $\B_{f_a}$, $\B_f$ with respect to the action of $\dd_t$, which was introduced in \cite{mic}. The {\it microlocal $V$-filtration} on $\BB_{f_a}$, $\BB_f$ is naturally associated with this, and the induced filtration on the structure sheaf $\OO_Y$, $\OO_{Y_a}$ gives the {\it microlocal multiplier ideals} $\Jt(\al X)$, $\Jt(\al X_a)$, which coincide with the usual multiplier ideals $\J(\al X)$, $\J(\al X_a)$ for $\al<1$, see (2.1) below. We then get the Thom-Sebastiani theorem for microlocal multiplier ideals which holds for any $\al\in\Q$ (see Theorem~(2.2) below), and this implies Theorem~1 by restricting to $\al<1$.
\sk
As another application of Theorem~2, we prove in \cite{MSS} a Thom-Sebastiani type theorem for the {\it spectral Hirzebruch-Milnor classes} introduced there. Such a result can be viewed as a global analogue of a similar assertion for the Steenbrink spectrum of hypersurface singularities, which follows from results in \cite{ScSt}, \cite{Va} in the isolated hypersurface singularity case, and from \cite{ts} or Theorem~2 in general.
(Finally, it is not clear to us how to relate the Thom-Sebastiani theorem in \cite{Ba} with the ones obtained in this paper.)
\sk
We thank M.~Musta\c{t}\v{a} for very important information concerning the relation between the summation formula and the Thom-Sebastiani type theorem for multiplier ideals, see Remark~(2.7) below. We also thank the referee for useful comments. The first named author is partially supported by NSF and NSA. The second named author is partially supported by Kakenhi 15K04816. The third named author is supported by the SFB 878 ``groups, geometry and actions''.
\sk
In Section 1 we first explain algebraic partial microlocalization together with microlocal $V$-filtration, and then prove Theorem~(1.2) below, which is a generalization of Theorem~2. In Section 2 we prove the Thom-Sebastiani type theorem for microlocal multiplier ideals in Theorem~(2.2) below.
\bs\bs
\vbox{\centerline{\bf 1. Thom-Sebastiani theorem}
\bsn
In this section we first explain algebraic partial microlocalization together with microlocal $V$-filtration, and then prove Theorem~(1.2) below, which is a generalization of Theorem~2.}
\msn
{\bf 1.1.~Algebraic microlocalization.} Let $Y$ be a smooth complex algebraic variety (or a connected complex manifold) with $f$ a non-constant function on $Y$, that is $f\in\Gamma(Y,\OO_Y)\setminus\C$. Let $i_f:Y\into Y\times\C$ be the graph embedding by $f$, and $t$ be the coordinate of the second factor of $Y\times\C$. Let $(M,F)$ be a holonomic filtered left $\D$-module (in particular, $\Gr^F_{\ssb}M$ is coherent over $\Gr^F_{\ssb}\D_Y$). Set
$$(M_f,F):=(i_f)_*^{\D}(M,F)=(M[\dd_t],F)\q\h{with}\q F_pM_f=\mopl_{i\in\Z}\,F_{p-i}M\sot\dd_t^i.
\leqno(1.1.1)$$
Here $(i_f)_*^{\D}$ denotes the direct image of filtered $\D$-modules, and the filtration $F$ is {\it indexed like a right $\D$-module}, see for instance \cite[Section 1.2]{MSS}, \cite[Section B.3]{ypg}. (In this case the filtration $F$ on the de Rham functor $\DR$ must be defined also as in the right $\D$-module case, see for instance \cite[1.2.2]{MSS}.) The above second isomorphism is as filtered $\OO_Y[\dd_t]$-modules, and the sheaf-theoretic direct image $(i_f)_*$ is omitted to simplify the notation.
The actions of $t$ and $\dd_{y_i}$ with $y_i$ local coordinates of $Y$ are given by
$$\aligned t(m\1\dd_t^j)&=fm\1\dd_t^j-jm\1\dd_t^{j-1},\\ \dd_{y_i}(m\1\dd_t^j)&=(\dd_{y_i}m)\,\dd_t^j-(\dd_{y_i}f)m\1\dd_t^{j+1}\q(m\in M).\endaligned
\leqno(1.1.2)$$
Here $m\1\dd_t^j$ is an abbreviation of $m\sot\dd_t^j$, or more precisely, $m\sot\dd_t^j\delta(t-f)$ with
$$\delta(t-f):=\h{$\frac{1}{t-f}$}\in\OO_{Y\times\C}\bl[\h{$\frac{1}{t-f}$}\br]/\OO_{Y\times\C}.$$
Note that $\delta(t-f)$ is also identified with $f^s$ (and $-\dd_tt$ with $s$), and this identification gives the relation between the Bernstein-Sato polynomial and the $V$-filtration, see for instance \cite{Ma2}.
\sk
Let $\M_f$ be the {\it algebraic partial microlocalization} of $M_f$ (see \cite{mic}), that is,
$$(\M_f,F)=(M[\dd_t,\dd_t^{-1}],F)\q\h{with}\q F_p\M_f=\mopl_{i\in\Z}\,F_{p-i}M\sot\dd_t^{i},
\leqno(1.1.3)$$
where the actions of $t$ and $\dd_{y_i}$ are given by (1.1.2).
\sk
In the case $M=\OO_Y$, we denote $M_f$, $\M_f$ respectively by $\B_f$, $\BB_f$, where the Hodge filtration $F$ on $\OO_Y$ is defined by $\Gr^F_p\OO_Y=0$ ($i\ne-d_Y$) with $d_Y:=\dim Y$ as in the case of right $\D$-modules, so that
$$F_p\B_f=\mopl_{0\les i\les p+d_Y}\,\OO_Y\sot\dd_t^i,\q F_p\BB_f=\mopl_{i\les p+d_Y}\,\OO_Y\sot\dd_t^i.
\leqno(1.1.4)$$
\sk
Let $V$ be the $V$-filtration of Kashiwara \cite{Ka2} and Malgrange \cite{Ma2} on $M_f$ along $t=0$, see for instance \cite[Proposition 1.9]{rat}. We assume that $V$ is indexed discretely by the {\it real numbers} $\R$ (instead of $\Q$ as in \cite{mhp}, \cite{mhm}), and moreover the following conditions are satisfied (as in the case of mixed Hodge modules, see \cite[Section 3.2.1]{mhp}):
$$\h{$\Gr^F_{\ssb}\Gr_V^{\al}M$ are coherent over $\Gr^F_{\ssb}\D_Y\q(\forall\,\al\in[0,1]$).}
\leqno(1.1.5)$$
$$t:V^{\al}(M_f,F)\simto V^{\al+1}(M_f,F)\q(\forall\,\al>0),
\leqno(1.1.6)$$
$$\dd_t:\Gr_V^{\al}(M_f,F)\simto\Gr_V^{\al-1}(M_f,F[-1])\q(\forall\,\al<1),
\leqno(1.1.7)$$
where $F[m]_p=F[m]^{-p}=F^{m-p}=F_{p-m}$ in general. We say that $(M,F)$ is $f$-{\it admissible} if conditions (1.1.6--7) hold, see \cite[Section 2.1]{ypg}.
\sk
Let $V$ be the microlocal $V$-filtration on $\M_f$ along $t=0$, see \cite{mic} in the case $M=\OO_Y$. In general this can be defined by modifying the $V$-filtration on $M_f$ so that
$$V^{\al}\M_f=\msum_{i\ges 0}\,\dd_t^{-i}V^{\al-i}M_f\q\q(\al\in\R),$$
where $M_f$ is naturally identified with a subsheaf of $\M_f$.
It is an exhaustive decreasing filtration on $\M_f$ indexed discretely by $\R$.
\sk
There is a canonical inclusion
$$\can:(M_f;F,V)\into(\M_f;F,V),$$
which is strictly compatible with $F$ by (1.1.1), (1.1.3). It is also strictly compatible with $V^{\al}$ for $\al<1$. Indeed, if there is $m\in M_f$ with $m\in V^{\al}\M_f$ for some $\al<1$, then there is some $k\gg 0$ with $\dd_t^km\in V^{\al-k}M_f$ by the above definition of $V$ on $\M_f$ (since $\dd_tV^{\al}M_f\subset V^{\al-1}M_f$), where $\,\can\,$ is omitted to simplify the notation. Then $m\in V^{\al}M_f$ by using the injectivity of (1.1.7) with $F$ forgotten. A similar argument implies the isomorphism (1.1.9) below with filtration $F$ forgotten, see also the proof of (1.1.9) below in the filtered case.
\sk
The filtration $V$ on $\M_f$ satisfies the following three conditions, and moreover it is uniquely determined by these.
\msn
(a) \,\,The $V^{\al}\M_f$ are locally finitely generated over $\D_Y[\dd_t^{-1}]$ ($\forall\,\al\in\R$).
\msn
(b) $\,\,t(V^{\al}\M_f)\subset V^{\al+1}\M_f,\,\,\,\dd_t(V^{\al}\M_f)=V^{\al-1}\M_f\,\,\,(\forall\,\al\in\R)$.
\msn
(c) \,\,The action of $\dd_tt-\al$ on $\Gr_V^{\al}\M_f$ is nilpotent ($\forall\,\al\in\R$).
\msn
The property~(a) follows (1.1.6) together with Nakayama's lemma applied to the $F_pV^{\al}M_f$. Indeed, these imply that the $V^{\al}M_f$ are locally finitely generated over $\OO_{Y\times\C}\langle\dd_{y_1},\dots,\dd_{y_{d_Y}}\rangle$, and any local section of $M_f$ is annihilated by a sufficiently high power of $t-f$.
\sk
By condition~(b) together with (1.1.3) there are canonical bi-filtered isomorphisms
$$\dd_t^{\,k}:(\M_f;F,V)\simto(\M_f;F[-k],V[-k])\q\q(\forall\,k\in\Z).
\leqno(1.1.8)$$
We then get the canonical filtered isomorphisms
$$\Gr_V^{\al}\can:\Gr_V^{\al}(M_f,F)\simto\Gr_V^{\al}(\M_f,F)\q\q(\forall\,\al<1).
\leqno(1.1.9)$$
Here it is enough to prove the strict surjectivity, since the injectivity with $F$ forgotten follows from the strict compatibility of the morphism $\,\can\,$ with $V^{\al}$ ($\al<1$) which is shown above. For any $m\in F_pV^{\al}\M_f$, we have for some $k\gg 0$ by the definition of $V$ on $\M_f$
$$\dd_t^km\in\can(V^{\al-k}M_f)\cap F_{p+k}\M_f=\can(F_{p+k}V^{\al-k}M_f),$$
where the last equality follows from the strictness of $\,\can:(M_f,F)\into(\M_f,F)$. Then (1.1.9) follows from (1.1.7--8).
\sk
For a holonomic filtered $\D_Y$-module $(M,F)$, we define the vanishing cycle functor by
$$\varphi_f(M,F):=\mopl_{\al\in(-1,0]}\,\Gr_V^{\al}(M_f,F)=\mopl_{\al\in(-1,0]}\,\Gr_V^{\al}(\M_f,F),
\leqno(1.1.10)$$
assuming (1.1.7) for the last isomorphism (using (1.1.9)), see Remark~(1.6)(i) below for the justification of the use of $\varphi_f$.
Here the filtration $F$ is {\it not} shifted, although it is shifted by $1$ using (1.1.8) if we consider $\Gr_V^{\al}\M_f$ for $\al\in(0,1]$ instead of $\al\in(-1,0]$. This definition is compatible with the one in \cite[5.1.3.3]{mhp} in the Hodge module case. Indeed, for $\al=0$, it coincides with the original definition of $\varphi_{f,1}$ in {\it loc.~cit.}, and we use (1.1.8--9) for the case $\al\in(-1,0)$, since $\psi_{f,\la}=\varphi_{f,\la}$ for $\la\ne 1$.
(In this section, we index the filtration $F$ like {\it right} $\D$-modules as is explained before (1.1.2), and $V$ is indexed decreasingly so that $V_{\al}=V^{-\al}$ and $\Gr^V_{\al}=\Gr_V^{-\al}$.)
\sk
We have a Thom-Sebastiani theorem as below. This is mentioned in \cite[Remark 4.5]{mic}, and follows from \cite{ts} in the Hodge module case (see \cite{DL}, \cite{GLM} for the motivic version, \cite{Ma1}, \cite{ScSt}, \cite{Va} for the isolated hypersurface singularity case with constant coefficients, and also \cite{Mas}, \cite{Sch1} for complexes with constructible cohomology sheaves). We give here a proof using algebraic partial microlocalization.
\msn
{\bf Theorem~1.2.} {\it Let $Y_a$ be a smooth complex algebraic variety $($or a connected complex manifold$)$ with $f_a$ a non-constant function, that is, $f_a\in\Gamma(Y_a,\OO_{Y_a})\setminus\C$, for $a=1,2$. Set $Y=Y_1\times Y_2$, and $f=f_1+f_2$ in the notation of the introduction. Let $(M_a,F)\,\,\,(a=1,2)$ be a holonomic filtered left $\D_{Y_a}$-module satisfying the assumptions {\rm(1.1.5--7)}, and define $\M_{a,f_a}$, $\M_f$ as in $(1.1)$ with $M:=M_1\boxtimes M_2$. Then there are canonical isomorphisms of holonomic filtered $\D_Y$-modules for $\al\in(-1,0]:$
$$\aligned\Gr^{\al}_V(\M_f,F)&=\mopl_{\al_1\in I(\al)}\,\Gr^{\al_1}_V(\M_{1,f_1},F)\boxtimes\Gr^{\al-\al_1}_V(\M_{2,f_2},F)\\
&\q\oplus\,\mopl_{\al_1\in J(\al)}\,\Gr^{\al_1}_V(\M_{1,f_1},F)\boxtimes\Gr^{\al-1-\al_1}_V(\M_{2,f_2},F[-1]),\endaligned
\leqno(1.2.1)$$
by replacing $Y_a$ with an open neighborhood of $X_a:=f_a^{-1}(0)$ in $Y_a$ {$(a=1,2)$} if necessary, where}
$$I(\al):=(-1,0]\cap[\al,\al+1)\q\q J(\al):=(-1,0]\cap[\al-1,\al).$$
\msn
{\bf Note.} Setting $\al_2=\al-\al_1$, $\,\al'_2=\al-1-\al_1$, we have
$$\al_1\in I(\al)\iff\al_1,\,\al_2\in(-1,0],\q\q\al_1\in J(\al)\iff\al_1,\,\al'_2\in(-1,0].
\leqno(1.2.2)$$
\msn
{\it Proof.} Set $\Si_a=\Sing\,f_a^{-1}(0)$ ($a=1,2$). By replacing $Y_a$ with an open neighborhood of $X_a$ ($a=1,2$) if necessary, we may assume
$$\Sing\,f^{-1}(0)=\Si_1\times\Si_2.$$
Indeed, $f^{-1}(0)$ is the inverse image of the anti-diagonal of $\C\times\C$ by $f_1\times f_2$.
\sk
We next show that there is the bi-filtered short exact sequence
$$0\to(\M_{1,f_1}\boxtimes\M_{2,f_2};F[1],V[1])\buildrel{\iota}\over\to(\M_{1,f_1}\boxtimes\M_{2,f_2};F,V)\buildrel{\eta}\over\to(\M_f;F,V)\to 0,
\leqno(1.2.3)$$
where $\iota$ is defined by
$$\dd_{t_1}\boxtimes id-id\boxtimes\dd_{t_2},$$
and $\eta$ by
$$\eta\bl(m_1\dd_{t_1}^{i_1}\boxtimes m_2\dd_{t_2}^{i_1}\br):=
(m_1\boxtimes m_2)\dd_t^{i_1+i_2}\q\h{for}\,\,\,\,m_a\in M_a\,\,(a=1,2),$$
(see also \cite[Section 4.1]{mic} for the case $M_a=\OO_{Y_a}$ ($a=1,2)$).
We have
$$(\M_{1,f_1}\boxtimes\M_{2,f_2};F,V):=(\M_{1,f_1};F,V)\boxtimes(\M_{2,f_2};F,V),\leqno(1.2.4)$$
where the external product $\boxtimes$ is taken as that of $\OO_{Y_a}$-modules ($a=1,2$). This implies the following filtered isomorphisms for $\al\in\R$ (see Lemma~(1.4) below for the proof):
$$\Gr_V^{\al}(\M_{1,f_1}\boxtimes\M_{2,f_2},F)\buildrel{\sim}\over\longleftarrow\mopl_{\al_1\in\R}\,\Gr_V^{\al_1}(\M_{1,f_1},F)\boxtimes\Gr_V^{\al-\al_1}(\M_{2,f_2},F).
\leqno(1.2.5)$$
By the definition of $\iota$ and by using (1.1.8) for $j=1$, the filtered isomorphism (1.2.5) implies that $\Gr_V^{\al}\iota$ is strictly injective for $F$ and we have the filtered isomorphism for $\al\in\R:$
$$({\rm Coker}\,\Gr_V^{\al}\iota,F)\cong\mopl_{\al_1\in(-1,0]}\,\Gr_V^{\al_1}(\M_{1,f_1},F)\boxtimes\Gr_V^{\al-\al_1}(\M_{2,f_2},F),
\leqno(1.2.6)$$
where the left-hand side is defined to be a quotient of $\Gr_V^{\al}(\M_{1,f_1}\boxtimes\M_{2,f_2},F)$.
\sk
We now show that the morphism $\eta$ in (1.2.3) induces an isomorphism of bi-filtered $\D_Y$-modules
$$({\rm Coker}\,\iota;F,V)\simto(\M_f;F,V),
\leqno(1.2.7)$$
which is also compatible with the action of $t$, $\dd_t$. Here the action of $t$ and $\dd_t$ on ${\rm Coker}\,\iota$ is defined respectively by $t_1+t_2$ and either $\dd_{t_1}$ or $\dd_{t_2}$ by the definition of $\iota$.
(To keep the exposition more concise, we choose not to formulate (1.2.7) as a theorem.)
The compatibility of the isomorphism (1.2.7) with $F$ follows from the definition (1.1.3). The compatibility with the filtration $V$ is equivalent to that $\eta$ is strictly compatible with the filtration $V$. By the uniqueness of the microlocal filtration $V$ explained in (1.1), this is also equivalent to that the quotient filtration $V$ on ${\rm Coker}\,\iota$ satisfies the conditions of the microlocal $V$-filtration in (1.1).
Here the finiteness condition (a) for $\M_f$ follows from that for $\M_{a,f_a}$.
Condition (b) follows from the definition of the action of $t$, $\dd_t$ explained above.
Condition (c) is verified also by using the definition of the action of $t$, $\dd_t$ on the left-hand side (especially $t=t_1+t_2$).
Thus (1.2.7) follows.
\sk
As is shown after (1.2.5), $\Gr_V^{\al}\iota$ is strictly injective for $F$. Hence $\iota$ is bistrictly injective for $F,V$. By (1.2.7), this implies that (1.2.3) is bistrictly exact for $F,V$, applying \cite[Theorem~1.2.12]{mhp} to
$$(\M_{1,f_1}\boxtimes\M_{2,f_2};F,V,G),$$
where the third filtration $G$ is defined so that $\Gr_i^G=0$ for $i\ne 0,1$ and $\Gr_0^G$ is given by the image of $\iota$. In particular, the cokernel commutes with $\Gr_V^{\al}$ in a compatible way with the filtration $F$. (This does not necessarily hold if $\iota$ is not bistrictly compatible with $F,V$.)
\sk
The assertion (1.2.1) now follows from (1.2.6--7). Indeed, (1.2.6) says that ${\rm Coker}\,\Gr_V^{\al}\iota$ for $\al\in(-1,0]$ is given by the direct sum over the index set defined by the conditions:
$$\al_1\in(-1,0],\q\al_1+\al_2=\al\in(-1,0]\q\h{with}\q\al_2:=\al-\al_1,$$
where $\al_2\in(-1,1)$, and does not necessarily belong to $(-1,0]$.
However, the difference with the union of the index sets of the direct sums in (1.2.1) (see also (1.2.2)) can be recovered by using (1.1.8) for $a=2$, $j=1$, where we get the shift of $F$ by $-1$ in the last term of (1.2.1).
Thus Theorem~(1.2) follows.
\msn
{\bf Remarks~1.3.} In the proof of $(1.2.5)$ given in Lemma~(1.4) below, we need the notion of convolution of filtrations and the related compatibility properties, which we now recall. In general, assume there are compatible $m$ filtrations $F_{(1)},\dots,F_{(m)}$ of an object $A$ of an abelian category $\A$ in the sense of \cite[Section 1]{mhp}. Here we assume that the inductive limit over a directed set is always an exact functor (for instance, the category of $\C$-vector spaces). Then we can show by using \cite[Theorem 1.2.12]{mhp} that the $m+1$ filtrations $F_{(1,2)},F_{(1)},\dots,F_{(m)}$ also form {\it compatible filtrations} of $A$, where $F_{(1,2)}$ is the {\it convolution} of $F_{(1)}$ and $F_{(2)}$, that is,
$$F_{(1,2)}^pA=\msum_{q\in\Z}\,F_{(1)}^qA\cap F_{(2)}^{p-q}A.$$
This assertion can be reduced to the finite sum case by using an inductive limit argument as above, and then to the finite filtration case (by replacing $F_{(1)}^p$, $F_{(2)}^q$ with $0$ for $p,q\gg 0$). Here we have the canonical isomorphisms
$$\Gr_{F_{(1,2)}}^pA\buildrel{\sim}\over\longleftarrow\mopl_{q\in\Z}\,\Gr_{F_{(1)}}^q\Gr_{F_{(2)}}^{p-q}A,$$
which is compatible with the filtrations $F_{(i)}$ ($i>2$).
This can be shown by using
$$\Gr_{F_{(2)}}^{p-q}\Gr_{F_{(1,2)}}^pA=\Gr_{F_{(1,2)}}^p\Gr_{F_{(2)}}^{p-q}A=\Gr_{F_{(1)}}^q\Gr_{F_{(2)}}^{p-q}A.$$
Here, taking an abelian category containing the exact category of $(m-2)$-filtered objects of $\A$ as in \cite[Sections 1.3.2--3]{mhp}, the assertion is essentially reduced to the case $m=2$, where it is more or less well-known.
\msn
{\bf Lemma~1.4.} {\it With the above notation and assumptions, we have the filtered isomorphism $(1.2.5)$.}
\msn
{\it Proof.} Let $F^{(a)}$ be the filtration on $\M_{1,f_1}\boxtimes\,\M_{2,f_2}$ induced by $F$ on $\M_{a,f_a}$, and similarly for $V_{(a)}$ ($a=1,2$). The filtrations $F$, $F^{(1)}$, $F^{(2)}$, $V$, $V_{(1)}$, $V_{(2)}$ on $\M_{1,f_1}\boxtimes\,\M_{2,f_2}$ form {\it compatible filtrations} in the sense of \cite[Section 1]{mhp}.
Note that $F$ is the {\it convolution} of $F^{(1)}$, $F^{(2)}$, and similarly for $V$, see Remark~(1.3) above for convolution.
We can prove the {\it compatibility} of the above six filtrations by using \cite[Theorem 1.2.12]{mhp}.
Indeed, the compatibility of the four filtrations $F^{(1)}$, $F^{(2)}$, $V_{(1)}$, $V_{(2)}$ follows from the definition, since the external product is an exact functor for both factors. Then we can apply Remark~(1.3) above, and (1.2.5) follows by taking $\Gr^{V_{(2)}}$, since we have the canonical isomorphism
$$\Gr_{V_{(1)}}^{\al_1}\Gr_{V_{(2)}}^{\al_2}(\M_{1,f_1}\boxtimes\M_{2,f_2})=\Gr_V^{\al_1}\M_{1,f_1}\boxtimes\Gr_V^{\al_2}\M_{2,f_2},$$
which is compatible with the filtrations $F^{(1)}$, $F^{(2)}$.
\sk
In the case of (1.2.5), however, there is an additional difficulty, since the filtration $V$ does not satisfy the condition $V^{\al}=0$ for $\al\gg 0$ (and similarly for $V_{(1)}^{\al}$, $V_{(2)}^{\al}$).
In order to avoid this problem, we can restrict to $V_{(1)}^{\beta}$, and take the inductive limit for $\beta\to-\infty$.
Note that the induced filtration $V_{(2)}$ on $\Gr_V^{\al}V^{\beta}_{(1)}$ satisfies the above property (since $V^{\,\gamma}_{(2)}\Gr_V^{\al}V^{\beta}_{(1)}=0$ for $\gamma>\al-\beta$). This finishes the proof of Lemma~(1.4).
\msn
{\bf 1.5.~Case $M=\OO_Y$.} From now on, we restrict to the case $M=\OO_Y$, and explain some well-known properties needed for applications to multiplier ideals.
\sk
Recall first that the morphism $\,\can\,$ in (1.1.9) {\it for\,} $\al=1$ is strictly surjective by \cite[Lemma 5.1.4 and Proposition 5.1.14]{mhp}. Indeed, by setting
$$X:=f^{-1}(0)\subset Y,$$
the morphism $\,\can\,$ in (1.1.9) for $\al=1$ is identified with the underlying morphism of filtered $\D_Y$-modules of the morphism $\,\can\,$ in the short exact sequence of mixed Hodge modules
$$0\to\Q_{h,X}[d_X]\to\psi_{f,1}\Q_{h,Y}[d_X]\buildrel{\rm can\,\,}\over\longrightarrow\varphi_{f,1}\Q_{h,Y}[d_X]\to 0.
\leqno(1.5.1)$$
(This follows for instance from \cite[Corollary 2.24]{mhm}.) Here
$$\Q_{h,X}:=a_X^*\Q\in D^b{\rm MHM}(X),$$
with $a_X:X\to pt$ the structure morphism (similarly for $\Q_{h,Y}$), and $\psi_{f,1}$, $\varphi_{f,1}$ are the unipotent monodromy part of the nearby and vanishing cycles functors \cite{De}. More generally, $\psi_{f,\la}$, $\varphi_{f,\la}$ for $\la\in\C$ can be defined for the underlying $\C$-complexes by using ${\rm Ker}(T_s-\la)$ with $T_s$ the semisimple part of the monodromy. (Note that the nearby and vanishing cycle functors preserve mixed Hodge modules {\it up to the shift of complex by} $1$ in this paper.)
\sk
Using (1.1.8--9) we have the canonical isomorphisms for any $\al\in\Q$
$$\DR_Y(\Gr_V^{\al}\BB_f)=\varphi_{f,\e(-\al)}\C_Y[d_X],
\leqno(1.5.2)$$
where $\e(-\al):=\exp(-2\pi i\al)$, and the left-hand side is the de Rham complex associated with the $\D_Y$-module $\Gr_V^{\al}\BB_f$ (viewed as an $\OO_Y$-module with an integrable connection) and shifted by $d_Y$ to the left; see also Remark~(1.6)(i) below for a generalization of (1.5.2).
\sk
The above argument implies that the morphism $\,\can:\Gr^1_V\B_f\to\Gr_V^1\BB_f$ is identified with the underlying morphism of the canonical morphism $\,\can\,$ in the short exact sequence (1.5.1), and hence the short exact sequence (2.1.9) below is identified with the short exact sequence deduced from (1.5.1) by restricting to $F_{-d_X}$ of the underlying filtered $\D_Y$-modules. 
\msn
{\bf Remarks 1.6.} (i) The isomorphism (1.5.2) holds for any $\al\in\R$ with $\BB_f$, $\C_Y[d_X]$ replaced respectively by $\M_f$, $\DR_Y(M)[-1]$ if $M$ is regular holonomic or if (1.1.5--7) are satisfied.
Indeed, this is well known for $M_f$ (instead of $\M_f$) if $\al\in[0,1)$ (see \cite[Proposition 3.4.12]{mhp} for the case (1.1.5--7) hold). Then it holds for $\M_f$ and for any $\al\in\R$ by (1.1.8--9).
\sk
(ii) In general the mixed Hodge modules are stable by $\psi[-1]$, $\varphi[-1]$. As a consequence, we get the shift of complex by $-1$ in the Thom-Sebastiani isomorphism (see also \cite{Mas}, \cite[Corollary 1.3.4 on p.~72]{Sch1}):
$$\varphi_f\C_Y=\varphi_{f_1}\C_{Y_1}\boxtimes\varphi_{f_2}\C_{Y_2}[-1],
\leqno(1.6.1)$$
which is closely related to the sign appearing in \cite[Theorems 3 and 4]{MSS}.
\sk
(iii) The Thom-Sebastiani isomorphism obtained by Theorem~(1.2) can be compatible with the one for complexes with constructible cohomology sheaves in \cite{Mas}, \cite{Sch1} via the de Rham functor in the mixed Hodge module case only after applying some automorphism of the external product. Indeed, it is expected that the former coincides with the one for the underlying filtered $\D$-module part in \cite{ts} (where the argument is not completely trivial), and it is shown in \cite[Proposition A.2]{Sch2} that the latter coincides with the one in \cite[Section 2]{ts}. However, these two Thom-Sebastiani isomorphisms in \cite{ts} can coincide only after applying an automorphism of the external product defined by using the beta function and the logarithm of the unipotent part of the monodromy as is seen in the definition of twisted external products, see {\it loc.~cit.}
\sk
(iv) By (1.1.4) we get
$$\Gr^F_p\Gr_V^{\al}\B_f=\Gr^F_p\Gr_V^{\al}\BB_f=0\q\h{if}\q p<-d_Y,\,\,\al\in(-1,0].
\leqno(1.6.2)$$
For $\al=0$, this uses the strict surjectivity of (1.1.9) for $\al=1$.
(It is closely related to the strict negativity of the roots of $b$-functions, see \cite{Ka1}.) This is used in \cite[3.1.7]{MSS} and also in Section 2 below (for instance in (2.1.9)) implicitly. 
\bs\bs
\vbox{\centerline{\bf 2. Application to multiplier ideals}
\bsn
In this section we prove the Thom-Sebastiani type theorem for microlocal multiplier ideals in Theorem~(2.2) below.}
\msn
{\bf 2.1.~Multiplier ideals.} Let $Y$ be a smooth complex algebraic variety (or a connected complex manifold), and $f$ be a non-constant function on $Y$, that is, $f\in\Gamma(Y,\OO_Y)\setminus\C$. Set $X:=f^{-1}(0)\subset Y$. Let $\J(\al X)\subset\OO_Y$ be the {\it multiplier ideal} of $X$ with coefficient $\al\in\Q$ (or $\R$ more generally). It can be defined by the local integrability of
$$|g|^2/|f|^{2\al}\q\h{for}\,\,\,\,g\in\OO_Y.
\leqno(2.1.1)$$
(For general references to multiplier ideals, see \cite{La}, \cite{Na}.) By definition, the $\J(\al X)$ form a decreasing sequence of ideal sheaves of $\OO_Y$ indexed by $\R$ and satisfying
$$\J(\al X)=\OO_Y\q(\al\les 0),\q \J((\al+1)X)=f\J(\al X)\q(\al\ges 0).
\leqno(2.1.2)$$
It is quite well known (and is quite easy to verify) that the multiplier ideals can be defined also by using an embedded resolution of $X\subset Y$, since the integrability condition can be expressed by the multiplicities of the pull-backs of $f,g$ and the differential form ${\rm d}y_1{\wedge}\cdots{\wedge}{\rm d}y_n$ (with $n=d_Y$) along the irreducible components of the exceptional divisor, see also \cite[1.1.1]{BS} for instance. This implies the {\it coherence} of the $\J(\al X)$ together with
$$\J(\al X)=\J((\al+\ep)X)\q(0<\forall\,\ep\ll 1).
\leqno(2.1.3)$$
The latter means that $\J(\al X)$ is {\it right-continuous} for $\al$. More precisely, for any $\al'\in\R$, the argument using an embedded resolution implies that we have for some $\beta,\beta'\in\Q$
$$\bl\{\al\in\R\mid\J(\al X)=\J(\al'X)\br\}{}=[\beta,\beta')\,\,\,\h{or}\,\,\,(-\infty,\beta').$$
\sk
We define the graded quotients $\G(\al X)$ by
$$\G(\al X):=\J((\al-\ep)X)/\J(\al X)\q(0<\ep\ll 1),$$
where the range of $\ep$ may depend on $\al$ (this is the same in (2.1.3)). We then have
$$\JC(X):=\bl\{\al\in\R\mid\G(\al X)\ne 0\br\}{}\subset\Q.$$
The members of $\JC(X)$ are called the {\it jumping coefficients} of $X$. We will restrict to rational numbers $\al$ when we consider $\J(\al X)$, $\G(\al X)$.
\sk
By (2.1.2) we get the isomorphisms
$$f:\G(\al X)\simto\G((\al+1)X)\q(\al>0),
\leqno(2.1.4)$$
and
$$\JC(X)=\bl(\JC(X)\cap(0,1]{\br)}\,+\,\N.
\leqno(2.1.5)$$
\sk
Consider the filtration $V$ on $\OO_Y$ induced by the filtration $V$ on $\B_f$ via the inclusion
$$\OO_Y=F_{-d_Y}\B_f\into\B_f.$$
By \cite[Theorem 0.1]{BS} we have
$$\aligned\J(\al X)&=V^{\al}\OO_Y&\h{if}\,\,\,\al\notin\JC(X),\\ \G(\al X)&=\Gr_V^{\al}\OO_Y=V^{\al}\OO_Y/\J(\al X)&\h{if}\,\,\,\al\in\JC(X).\endaligned
\leqno(2.1.6)$$
The last isomorphism is a consequence of the assertion that $\J(\al X)$ is {\it right-continuous} for $\al$ as is explained above, although $V^{\al}\OO_Y$ is {\it left-continuous} for $\al$.
\sk
We now consider the {\it microlocal} $V$-filtration on $\OO_Y$ which is denoted by $\Vt$, and is induced by the filtration $V$ on $\BB_f$ via the isomorphism
$$\OO_Y=\Gr^F_{-d_Y}\BB_f.$$
\vskip-2mm\nin
Set
$$\JCt(X):=\bl\{\al\in\Q\mid\Gr_{\Vt}^{\al}\OO_Y\ne 0\br\}.$$
We have by (1.1.9)
$$\JCt(X)\subset(0,+\infty),\q\JCt(X)\cap(0,1)=\JC(X)\cap(0,1).
\leqno(2.1.7)$$
However, the last equality does not necessarily hold if $(0,1)$ is replaced by $(0,1]$ (since $\JCt(X)$ does not necessarily contain $1$), and (2.1.5) with $\JC(X)$ replaced by $\JCt(X)$ does not necessarily hold, see Example~(2.6)(ii) below.
\sk
We have the {\it microlocal multiplier ideals} $\Jt(\al X)$ together with their graded quotients $\Gt(\al X)$ such that $\Jt(\al X)$ is {\it right-continuous} (see a remark after (2.1.3)) and
$$\aligned\Jt(\al X)&={\Vt}^{\al}\OO_Y&\h{if}\,\,\,\al\notin\JCt(X),\\ \Gt(\al X)&=\Gr_{\Vt}^{\al}\OO_Y={\Vt}^{\al}\OO_Y/\Jt(\al X)&\h{if}\,\,\,\al\in\JCt(X).\endaligned
\leqno(2.1.8)$$
As for the relation with the usual multiplier ideals, we have the following short exact sequence in the notation of (1.1):
$$0\to\omt_X{\otimes_{\OO_X}}\om_X^{\vee}\to F_{-d_X}\Gr_V^1\B_f\buildrel\can\over\longrightarrow F_{-d_X}\Gr_V^1\BB_f\to 0,
\leqno(2.1.9)$$
where $X$ is assumed to be reduced, and
$$\omt_X:=(\rho)_*\om_{\X}\into\om_X,\q\om_X^{\vee}:=\Hc om_{\OO_X}(\om_X,\OO_X),$$
with $\rho:\X\to X$ a resolution of singularities (see \cite[4.1.1]{MSS}). We have
$$\om_X=\om_Y\otimes_{\OO_Y}\OO_X,$$
since $X$ is globally defined by $f$, see for instance \cite[Lemma 2.9]{rat}. The coherent sheaf $\om_X/\omt_X$ may be called the ``irrationality" of the singularities of $X$, see also \cite[4.2.5]{MSS}.
\sk
By (1.1.9), (2.1.9) we then get
$$\J(\al X)=\Jt(\al X),\q\G(\al X)=\Gt(\al X)\q(\al<1),
\leqno(2.1.10)$$
$$\Jt(X)/\J(X)=\omt_X{\otimes_{\OO_X}}\om_X^{\vee}\subset\OO_X\q(\al=1),
\leqno(2.1.11)$$
$$0\to\omt_X{\otimes_{\OO_X}}\om_X^{\vee}\to\G(X)\to\Gt(X)\to 0\q(\al=1).
\leqno(2.1.12)$$
Here we assume $X$ reduced in (2.1.11--12).
Note that we have by (2.1.2)
$$\J(X)=\OO_Y(-X)=\I_X\q(\al=1),
\leqno(2.1.13)$$
where the last term is the ideal sheaf of $X$.
\ms
We have the Thom-Sebastiani type theorem for microlocal multiplier ideals as follows.
\msn
{\bf Theorem~2.2.} {\it With the notation and assumption of Theorem~$(1.2)$, there are equalities for any $\al\in\Q:$
$$\Jt(\al X)=\msum_{\al_1+\al_2=\al}\,\Jt(\al_1X_1)\boxtimes\Jt(\al_2X_2)\q\h{in}\,\,\,\,\OO_Y=\OO_{Y_1}\boxtimes\OO_{Y_2}.
\leqno(2.2.1)$$
by replacing $Y_a$ with an open neighborhood of $X_a=f_a^{-1}(0)$ in $Y_a$ {$(a=1,2)$} so that $\Si=\Si_1\times\Si_2$ if necessary. Here we may assume $\al_1,\al_2\in(0,\al)$ by the first equality in $(2.1.2)$ together with $(2.1.10)$.}
\msn
{\it Proof.} In (2.2.1) we may replace $\al_1+\al_2=\al$ by $\al_1+\al_2\ges\al$, and assume for $0<\ep\ll 1/m$
$$\al_a\in\JCt(X_a)-\ep\,\,\,(a=1,2),
\leqno(2.2.2)$$
(since $\Jt(\al X)$ is right-continuous), where $m$ is a positive integer such that $\JCt(X_a)\in\Z/m$.
We now show that (2.2.1) is equivalent to the following.
$$\Vt^{\al}\OO_Y=\msum_{\al_1+\al_2=\al}\,\Vt^{\al_1}\OO_{Y_1}\boxtimes\Vt^{\al_2}\OO_{Y_2}\q\h{in}\,\,\,\,\OO_Y=\OO_{Y_1}\boxtimes\OO_{Y_2}.
\leqno(2.2.3)$$
We may replace $\al_1+\al_2=\al$ by $\al_1+\al_2\ges\al$ in (2.2.3), and assume
$$\al_a\in\JCt(X_a)\,\,\,(a=1,2),
\leqno(2.2.4)$$
since $\Vt^{\al}$ is left-continuous. However, we may also assume (2.2.2) with $0<\ep\ll 1/m$ instead of (2.2.4) by replacing $\al$ with $\al-2\ep$ if necessary. (Here $\ep$ may depend on $\al$.) The equivalence between (2.2.1) and (2.2.3) then follows from (2.1.8).
\sk
We can show (2.2.3) by taking $\Gr_{-d_Y}^F$ of the isomorphism (1.2.7) and calculating $\Gr^F$ of $\iota$ in (1.2.7), since $\Gr^F_{p_a}(\BB_{f_a},V)$ is essentially independent of $p_a$ by (1.1.8) ($a=1,2$).
This finishes the proof of Theorem~(2.2).
\msn
{\bf Corollary 2.3.} {\it With the notation and assumption of Theorem~$(1.2)$, there are canonical isomorphisms for any $\al\in\Q:$
$$\Gt(\al X)=\mopl_{\al_1+\al_2=\al}\,\Gt(\al_1X_1)\boxtimes\Gt(\al_2X_2),
\leqno(2.3.1)$$
by replacing $Y_a$ with an open neighborhood of $X_a=f_a^{-1}(0)$ in $Y_a$ {$(a=1,2)$} if necessary. Here we may assume $\al_1,\al_2\in(0,\al)$ as in Theorem~$(2.2)$.}
\ms
It is also possible to deduce this from Theorem~(1.2).
Combining Corollary~(2.3) with (2.1.10), (2.1.12), we get the following.
\msn
{\bf Corollary~2.4.} {\it With the notation and assumption of Theorem~$(1.2)$, assume further $X$ reduced. We have the short exact sequence for $\al=1:$
$$0\to\omt_X{\otimes_{\OO_X}}\om_X^{\vee}\to\G(X)\to\mopl_{\al_1+\al_2=1}\,\G(\al_1X_1)\boxtimes\G(\al_2X_2)\to 0.
\leqno(2.4.1)$$
by replacing $Y_a$ with an open neighborhood of $X_a=f_a^{-1}(0)$ in $Y_a$ {$(a=1,2)$} if necessary.
Here we may assume $\al_1,\al_2\in(0,1)$ as in Theorem~$(2.2)$.}
\msn
{\bf 2.5.~Proof of Theorem~1.} The assertion follows from Theorem~(2.2) and Corollary~(2.3) together with (2.1.10).
\msn
{\bf Examples~2.6.} (i) Let $Y=\C$ with coordinate $z$. Set $f=z^m$ for $m\ges 2$. Then
$$\Vt^{i/m}\OO_Y=\OO_Yz^k\q\h{if}\q i=k+1+\bl[k/(m-1)\br].
\leqno(2.6.1)$$
Indeed, we have
$$V^{i/m}\OO_Y=\OO_Yz^{i-1}\q(i\in[1,m-1]),
\leqno(2.6.2)$$
where $V$ is the usual $V$-filtration on $\OO_Y$, see (2.1). This is compatible with \cite{exN}, and can be proved by using the multiplier ideals together with (2.1.6).
\sk
We then get the inclusion $\supset$ in (2.6.1) by using (2.6.2) together with the definition of the action of $\dd_z$ in (1.1.2) and (1.1.8), since $\Gr^F\dd_z$ preserves the filtration $\Vt$ and $\dd_zf=mz^{m-1}$. So it is enough to show
$$\dim \Gr_{\Vt}^{i/m}\OO_Y=\begin{cases}1&\h{if}\,\,\,i\ges 1,\,\,i/m\notin\Z,\\ 0&\h{otherwise}.\end{cases}
\leqno(2.6.3)$$
We have
$$\Ht^0(F_{\!f,0},\C)_{\la}=\begin{cases}\C&\h{if}\,\,\,\la^m=1\,\,\,\h{and}\,\,\,\la\ne 1,\\
\,0&\h{otherwise,}\end{cases}$$
where $\Ht$ is the reduced cohomology, and $F_{\!f,0}$ is the Milnor fiber so that $\varphi_{f,\la}\C_Y$ is identified with $\Ht^0(F_{\!f,0},\C)_{\la}$ in this case, see also (1.5.2). The assertion (2.6.3) then follows by using (1.1.8) and recalling the definition of the direct image of filtered $\D$-modules by the inclusion $\{0\}\into\C$, see for instance \cite[Section 1.2]{MSS}, \cite[Section B.3]{ypg}.
\sk
(ii) Let $Y=\C^d$ with coordinates $z_1,\dots,z_d$. Set $f=\sum_{j=1}^dz_j^{m_j}$ for $m_j\ges 2$ ($j\in[1,d]$). Then Example~(i) together with (2.2.3) implies
$$\Vt^{\al}\OO_Y=\msum_{\nu}\OO_Yz^{\nu},
\leqno(2.6.4)$$
where the summation is taken over $\nu=(\nu_1,\dots,\nu_d)\in\N^d$ satisfying
$$\h{$\msum_{j=1}^d\,\frac{1}{\,m_j}\bl(\nu_j+1+\bl[\frac{\nu_j}{m_j-1}\br]\br){}\ges\al$}.
\leqno(2.6.5)$$
In particular, we have
$$\Vt^{\alt_f}\OO_Y=\OO_Y\ne\Vt^{>\alt_f}\OO_Y\,\,\,\,\h{with}\,\,\,\,\alt_f:=\msum_{j=1}^d\,\h{$\frac{1}{\,m_j}$},\,\,\,\,\lct(X)=\min\{1,\alt_f\}.
\leqno(2.6.6)$$
\sk
By (2.6.4--5) we see that the microlocal $V$-filtration on $\OO_{Y,0}$ has nothing to do with the filtration $V$ on the {\it microlocal} Gauss-Manin system as in \cite{exN}. Indeed, the latter coincides with the usual Gauss-Manin system (since the Milnor fiber is contractible), and the filtration $V$ on it is induced by the {\it usual} $V$-filtration on $\B_f$, see \cite[Proposition 3.4.8]{mhp}.
\msn
{\bf Remark~2.7.} Theorem~1 and Corollary~1 were essentially known to specialists according to M.~Musta\c{t}\v{a}. (Here we use the multiplicative notation of multiplicative ideals following him.) Indeed, setting $\af:=(f_1,f_2)\subset\OO_Y$, the summation formula \cite[Theorem 0.3]{Mu} implies
$$\J(\af^{\al})=\msum_{\al_1+\al_2=\al}\,\J\bl((f_1)^{\al_1}(f_2)^{\al_2}\br),
\leqno(2.7.1)$$
and we have by \cite[Proposition 9.5.22]{La}
$$\J\bl((f_1)^{\al_1}(f_2)^{\al_2}\br)=\J\bl((f_1)^{\al_1}\br)\boxtimes\,\J\bl((f_2)^{\al_2}\br).
\leqno(2.7.2)$$
For each $\al\in(0,1)$, it follows from \cite[Proposition 9.2.28]{La} (with $k=1$) that we have for $c_1,c_2\in\C$ {\it general}
$$\J(\af^{\al})=\J\bl((c_1f_1+c_2f_2)^{\al}\br).
\leqno(2.7.3)$$
We thus get the following equality of ideal sheaves on $Y$ for $c_1,c_2\in\C$ general:
$$\J\bl((c_1f_1+c_2f_2)^{\al}\br)\,=\msum_{\al_1+\al_2=\al}\,\J\bl((f_1)^{\al_1}\br)\boxtimes\,\J\bl((f_2)^{\al_2}\br).
\leqno(2.7.4)$$
This may be viewed effectively as a Thom-Sebastiani type theorem for multiplier ideals. Here it is not necessarily immediate that the equality holds with $c_1=c_2=1$, that is, we have the equality 
$$\J\bl((c_1f_1+c_2f_2)^{\al}\br)=\J\bl((f_1+f_2)^{\al}\br).
\leqno(2.7.5)$$
This can be shown easily in the case when $f_1$ or $f_2$ is analytic-locally a weighted homogeneous polynomial by using the (local) $\C^*$-action on $Y_1$ or $Y_2$ (together with GAGA). Indeed, we have $\lambda^*f_1=(c_2/c_1)f_1$ with $\lambda^*\J\bl((f_1)^{\al_1}\br)=\J\bl((f_1)^{\al_1}\br)$ for some $\lambda\in\C^*$ if $f_1$ is a weighted homogeneous polynomial. In general, the assertion can be reduced to the above case by taking a resolution of $f_1$ or $f_2$, and applying \cite[Theorem 9.2.33]{La}.
\msn
{\bf Remark~2.8.} The {\it Hodge ideals\,} in \cite{MP} coincide with the microlocal multiplier ideals modulo the ideal of hypersurface, and hence the $j$-{\it log canonicity} ({\it loc.~cit.})\ is determined by $\min\JCt(X)$, the minimal microlocal jumping coefficient, which coincides with $\alt_f$, the maximal root of the microlocal $b$-function up to a sign, see \cite{him}.

\end{document}